\documentclass[leqno]{article}
\usepackage{graphicx} 
\usepackage{amsfonts}
\usepackage{amsmath}
\usepackage{amsthm}
\usepackage{amssymb}
\usepackage[symbol]{footmisc}
\usepackage{enumerate}
\usepackage{chngcntr}
\usepackage{sectsty}

\counterwithin*{equation}{section}
\sectionfont{\large}

\newcommand{\calA}{{\mathcal{A}}}

\newcommand{\calE}{{\mathcal{E}}}
\newcommand{\calF}{{\mathcal{F}}}

\newcommand{\calK}{{\mathcal{K}}}

\newcommand{\calM}{{\mathcal{M}}}

\newcommand{\calP}{{\mathcal{P}}}

\newcommand{\calS}{{\mathcal{S}}}

\newcommand{\calU}{{\mathcal{U}}}

\newcommand{\N}{{\mathbb{N}}} 
\newcommand{\R}{{\mathbb{R}}}

\title{On the measurability of a numerical function with respect to a family of sets}
\author{Gabriele H. Greco \thanks{Author's address: Department of Mathematics, University of Trento, 38050 Povo (Trento).} }
\date{}

\begin{document}

\noindent The following document is a translation (from French to English) of:

\bigskip \bigskip \bigskip

\noindent{\LARGE\textsc{Gabriele~H.~Greco}}

\bigskip

\noindent{\Large\bf Sur la mesurabilit\'e d'une fonction num\'erique par rapport \`a une famille d'ensembles} 

\bigskip

\noindent{\large \it Rendiconti del Seminario Matematico della Universit\`a di Padova}, tome 65 (1981), pp. 163--176.

\bigskip \bigskip \bigskip

\noindent Translated by: Jonathan M. Keith, School of Mathematics, Monash University, jonathan.keith@monash.edu

\bigskip

\noindent With thanks to: Prof. Andrea D'Agnolo, Editor-in-Chief of the above journal, for permission to publish this translation.

\newpage

\maketitle

\section*{Introduction.}

In this introduction and in the sequel, we will denote by the symbols $X$, $\calP(X)$, $[0, +\infty]^X$, $\overline{\R}^X$, $\N$ respectively, a non-empty set, the class of subsets of $X$, the class of functions defined on $X$ with values in $[0, +\infty]$, the class of functions defined on $X$ with values in the extended real line $\overline{\R}$, and the set of integers $n > 0$. Following the notations established by P. A. Meyer in~\cite{meyer1966}, we will call a {\em paving} on $X$ any family of sets $\calE \subset \calP(X)$ that contains the empty set; we will call a {\em paved set} the pair $(X, \calE)$ comprised of a set $X$ and a paving $\calE$ on $X$; the phrase ``$\calE$ is stable for $(\cup f, \cap d)$'' (and likewise analogous sentences) signifies that $\calE$ is stable for finite unions and countable intersections. We will also use the letters of the Greek alphabet $\alpha$, $\beta$, $\gamma$, $\delta$ to designate monotonic set functions, that is, functions defined on $\calP(X)$ with values in $[0, +\infty]$ such that $\alpha(\emptyset) = 0$ and $\alpha(A) \leqslant \alpha(B)$ if $A \subset B$. By $\int_X f d\alpha$ we will designate the real number, finite or not, called the {\em integral of $f$ over $X$ with respect to the monotone set function $\alpha$}, defined by $\int_X f d\alpha = \int_0^{+\infty} \alpha\{f > t\}dt$ if $f \in [0, +\infty]^X$, and by $\int_X f d\alpha = \int_X f^+ d\alpha - \int_X f^- d\alpha$ if $f \in \overline{\R}^X$ and if $\int_X f^+ d\alpha < \infty$ or $\int_X f^- d\alpha < \infty$; if $A \in \calP(X)$, we set $\int_A f d\alpha = \int_X f \varphi_A d\alpha$, where $\varphi_A$ is the characteristic function of the set $A$. These integrals have been studied independently in~\cite{degiorgi1977} and in~\cite{greco1977}; for their properties we will refer to~\cite{greco1977}. If $\calF$ is a filter of $\calP(X)$ we will denote by $\alpha_{\calF}$ and $\beta_{\calF} : \calP(X) \rightarrow \{0, 1\}$ the set functions defined by ``$\alpha_{\calF}(A) = 1$ if and only if $A \in \calF$'' and ``$\beta_{\calF}(A) = 1$ if and only if $A \cap F \neq \emptyset$ for any $F \in \calF$''; for any filter $\calF$ of $\calP(X)$ and for any $f \in [0, +\infty]^X$ it holds that $\int_X f d\alpha_{\calF} = \liminf_{\calF} f$\footnote{Translator's note: The original uses `minlim' for `liminf' and `maxlim' for `limsup' here and throughout.} and $\int_X f d\beta_{\calF} = \limsup_{\calF} f$ (see example $3^{\circ}$ in~\cite{greco1977}). We set $0 \cdot \infty = \infty \cdot 0 = 0$ and $\inf \emptyset = +\infty$. If $(X, \calE)$ is a paved set, we will denote by $\calS^+(X, \calE)$ the class of functions $\sum_{i = 1}^{\infty} a_i \varphi_{H_i}$, where $a_i \in [0, +\infty]$, $\{ H_i \}_{i \in \N} \subset \calE$ and $H_i \supset H_{i + 1}$ for all $i \in \N$.

In this article we propose and study a definition of measurability of a numerical function with respect to a family of sets $\calE \subset \calP(X)$ that meets this requirement: ``determine numerical functions $f$ such that the integral $\int_X f d\alpha$ depends only on those values of $\alpha$ on $\calE$, regardless of the monotone set function $\alpha: \calP(X) \rightarrow [0, +\infty]$''; that is, we say that a function $f \in [0, +\infty]^X$ is $\calE$-measurable if one always has $\int_X f d\alpha = \int_X f d\beta$, whenever the monotone set functions $\alpha, \beta$ are such that $\alpha(H) = \beta(H)$ for all $H \in \calE$. The function $f$ is then $\calE$-measurable if and only if for each pair $a, b \in (0, +\infty)$ with $a > b$, there exists a set $H \in \calE$ such that $\{f \geqslant a\} \subset H \subset \{f>b\}$ (see Th.~1). This definition permits each paved set $(X, \calE)$ to be associated with the class $\calM(X, \calE)$ of $\calE$-measurable functions, which can be: the well known class of measurable functions with respect to a $\sigma$-algebra (if $\calE$ is a $\sigma$-algebra), the class of functions without oscillatory discontinuities (if $X = \R$ and $\calE$ is the set algebra generated by the intervals of $\R$; see Remark~2), the class of continuous functions (if $X$ is a compact and totally disconnected Hausdorff topological space and if $\calE$ is the algebra of clopen sets), the class of functions which are continuous almost everywhere (if $X = [0, 1]$ and $\calE$ is the algebra of Peano-Jordan measurable sets), the class of $\alpha$-measurable sets in the sense of Carath\'eodory, with respect to a set function $\alpha$ (if $\calE$ is the algebra of $\alpha$-measurable sets in the sense of Carath\'eodory; see Remark~5); in the case where $\calE$ is an algebra, it follows that a function $f \in \overline{\R}^X$ is $\calE$-measurable if and only if $f$ extends by continuity to the Stone space associated with the algebra $\calE$.

In this article, we also notice analogies between pseudo-compact topological spaces and semi-compact paved sets (see s.~2), and between normal topological spaces and sets on which there are two good pavings (see s.~3). 

\section{$\calE$-measurable functions.}

{\it Let $(X, \calE)$ be a paved set. We say that a function $f \in [0, +\infty]^X$ is $\calE$-measurable if we always have $\int_X f d\alpha = \int_X f d\beta$, in the case where $\alpha$ and $\beta$ are monotonic set functions such that $\alpha(H) = \beta(H)$ for all $H \in \calE$. We also say that $f \in \overline{\R}^X$ is $\calE$-measurable if $f^+$ and $f^-$ are $\calE$-measurable. }

We will denote by $\calM(X, \calE)$ the class of functions $f \in \R^X$ which are $\calE$-measurable and we set $\calM^+(X, \calE) = \calM(X, \calE) \cap [0, +\infty]^X$.

The simplest $\calE$-measurable functions are the functions that belong to \linebreak $\calS^+(X, \calE)$; indeed if $\{ a_i \}_{i\in \N} \subset [0. +\infty]$; $\{ H_i \}_{i \in \N} \subset \calE$ and $H_{i + 1} \subset H_i$ for all $i \in \N$, we have $\int_X \left( \sum_{i=1}^{\infty} a_i \varphi_{H_i} \right)  d\alpha = \sum_{i=1}^{\infty} a_i \alpha(H_i)$ for any monotonic set function $\alpha$.

{\it \textsc{Lemma 1}. Let $(X, \calE)$ be a paved set and let $f \in [0, +\infty]^X$. The following conditions are equivalent:
\begin{enumerate}[i)]
\item $f$ is $\calE$-measurable,

\item for each pair $\alpha, \beta: \calP(X) \rightarrow \{0, 1\}$ of monotonic set functions such that $\alpha(H) = \beta(H)$ for all $H \in \calE$, we have
\[
\int_X f d\alpha = \int_X f d\beta.
\]
\end{enumerate}}

\textsc{Proof}. i) $\implies$ ii): The condition ii) is a particular case of condition i). ii) $\implies$ i): let us first set, for any monotonic set function $\alpha: \calP(X) \rightarrow [0, +\infty]$ and for any $s \in (0, \infty)$ and $H \in \calP(X)$, $\alpha_s(H) = 1$ or $0$ if we have respectively $\alpha(H) > s$ or $\alpha(H) \leqslant s$. If $\alpha, \beta: \calP(X) \rightarrow [0, +\infty]$ are two montonic set functions such that $\alpha(H) = \beta(H)$ for all $H \in \calE$, we also have $\alpha_s(H) = \beta_s(H)$ for all $s \in (0, +\infty)$. According to ii), we therefore have $\int_X f d\alpha_s = \int_X f d\beta_s$ for all $s \in (0, +\infty)$. We then deduce that
\[
\int_X f d\alpha = \lim_n \sum_{i=1}^{\infty} \frac{1}{2^n} \int_X f d\alpha_{i/2^n} = \lim_n \sum_{i=1}^{\infty} \frac{1}{2^n} \int_X f d\beta_{i/2^n} = \int_X f d\beta.
\]

{\it \textsc{Lemma 2}. Let $(X, \calE)$ be a paved set. If the function $f \in [0, +\infty]^X$ is $\calE$-measurable, then the functions $\alpha f$, $f \wedge a$, $f \vee a - a$ are also $\calE$-measurable for all $a \in [0, +\infty)$.}

\textsc{Proof}. For any monotonic set function $\alpha$ defined on $\calP(X)$ with values in $\{ 0, 1 \}$ it follows that:
\[
\int_X (f \wedge a) d\alpha = \left( \int_X f d\alpha \right) \wedge a \mbox{ and } \int_X (f \vee a - a) d\alpha = \left( \int_X f d\alpha \right) \vee a - a 
\]
for any $a \in [0, +\infty)$ and $f \in [0, +\infty]^X$. We therefore deduce Lemma~2 from Lemma~1. $\quad \quad \quad \blacksquare$

{\it \textsc{Theorem 1}. Let $(X, \calE)$ be a paved set and let $f \in [0, +\infty]^X$. The following conditions are equivalent:
\begin{enumerate}[i)]
\item $f$ is $\calE$-measurable,

\item for each pair $a, b \in (0, +\infty)$ such that $a > b$, there exists a set $H \in \calE$ such that $\{ f \geqslant a \} \subset H \subset \{ f > b \}$,

\item there exists a sequence $\{ g_n \}_{n \in \N} \subset \calS^+(X, \calE)$ such that $g_n \leqslant f$ for all $n \in \N$ and the sequence $\{ g_n \wedge a \}_{n \in \N}$ converges uniformly to $f \wedge a$ for all $a \in (0, +\infty)$.
\end{enumerate}}

\textsc{Proof}. i) $\implies$ ii): we reason by contradiction. If $f$ does not satisfy condition ii), there exist two real numbers $a, b \in (0, +\infty)$ with $a > b$ such that
\begin{equation}
\nexists H \in \calE \quad \mbox{such that } \{ f > a \} \subset H \subset \{ f > b \}\footnote{Translator's note: The original has $\{ f \geqslant a \} \subset H \subset \{ f > b \}$, but the proof seems to require this equivalent condition.}.
\end{equation}
Let us then define the two set functions $\tau_1, \tau_2: \calP(X) \rightarrow [0, 1]$ by
\[
\tau_1(A) = \left\{ 
\begin{array}{ll}
0 & \mbox{if } A = X \\
\sup \{ g(x): x \notin A \} & \mbox{if } A \neq X
\end{array}
\right.
\]
and
\[
\tau_2(A) = \left\{ 
\begin{array}{ll}
1 & \mbox{if } A = \emptyset \\
\inf \{ g(x): x \in A \} & \mbox{if } A \neq \emptyset,
\end{array}
\right.
\]
where $g = (f \wedge a - f \wedge b) / (a - b)$. According to this definition of $\tau_1, \tau_2$, in virtue of (1) it follows that\footnote{Translator's note: The original has $\tau_1(H) < t$ in~(3) and $\tau_2(H) > t$ in~(4).}
\begin{align}
& \tau_1(\emptyset) = \tau_2(\emptyset) = 1 \mbox{ and the functions } \tau_1, \tau_2 \mbox{ are decreasing;} \\
& \tau_1(H) \leqslant t \mbox{ and } \tau_2(H) = 0, \mbox{ if } H \in \calE, H \supset \{ g > t \} \mbox{ and } t \in (0, 1); \\
& \tau_1(H) = 1 \mbox{ and } \tau_2(H) \geqslant t, \mbox{ if } H \in \calE, H \subset \{ g > t \} \mbox{ and } t \in (0, 1).
\end{align}
Now the two monotonic set functions $\alpha, \beta: \calP(X) \rightarrow [0, +\infty]$, thus defined
\begin{align*}
\alpha(A) &= \sup \{ 2 - \tau_1(H) - \tau_2(H): H \in \calE \mbox{ and } H \subset A \}, \\
\beta(A) &= \inf \{ 2 - \tau_1(H) - \tau_2(H): H \in \calE \mbox{ and } H \supset A \},
\end{align*}
are therefore equal on $\calE$ and $\beta \geqslant \alpha$, according to (2). Moreover, we have $\beta \{ g > t \} > 1 > \alpha \{ g > t \}$ for all $t \in (0, 1)$, according to (3) and (4); this proves that $\int_X g d\beta > \int_X g d\alpha$, that is the function $g = (f \wedge a - f \wedge b) / (a - b)$ is not $\calE$-measurable, according to the definition of $\calE$-measurability. But this contradicts the hypothesis that $f$ is $\calE$-measurable; indeed from the $\calE$-measurability of $f$ we must deduce that $g$ must be $\calE$-measurable, according to Lemma~2. We thus complete the proof of i) $\implies$ ii).

ii) $\implies$ iii): if $f \in [0, +\infty]^X$ satisfies condition ii), for any $i, n \in \N$ there exists a set $H_{i, n} \in \calE$ such that 
\[
\left\{ f \geqslant \frac{i + 1}{2^n} \right\} \subset H_{i, n} \subset \left\{ f > \frac{i}{2^n} \right\} \quad \mbox{and} \quad H_{i + 1, n} \subset H_{i, n}
\]
We therefore have
\[
f \geqslant \frac{1}{2^n} \sum_{i = 1}^{\infty} \varphi_{H_{i,n}} \geqslant f \wedge \frac{1}{2^{n-1}} - \frac{1}{2^{n-1}}.
\]
The sequence $\{ g_n \}_{n \in \N}$, where $g_n = 1/2^n \sum_{i=1}^{\infty} \varphi_{H_{i,n}}$, satisfies condition iii).


iii) $\implies$ i): the $\calE$-measurability of all $g \in \calS^+(X, \calE)$ and the properties of the integral $\int_X \mbox{--- } d\alpha$ (see~\cite{greco1977}, Prop.~1) entail the $\calE$-measurability of $f$. $\quad \quad \quad \blacksquare$

\textsc{Remark 1}. Let $(X, \calE)$ be a paved set, let $\delta: \calE \rightarrow [0, +\infty]$ be a set function such that $\delta(\emptyset) = 0$ and $\delta(A) \leqslant \delta(B)$ for all $A, B \in \calE$ with $A \subset B$. The preceding theorem makes it possible to define in a satisfying manner the integral $\int_X \mbox{--- } d\delta$; it is defined by setting $\int_X f d\delta = \int_X f d\alpha$, where $\alpha: \calP(X) \rightarrow [0, +\infty]$ is a monotonic set function such that $\alpha(H) = \delta(H)$ for all $H \in \calE$, when $f \in \calM(X, \calE)$ and $\int_X f^+ d\alpha < +\infty$ or $\int_X f^- d\alpha < +\infty$. According to this definition of the integral we can analyse the properties of $\delta$ by means of the properties of $\int_X \mbox{--- } d\alpha$. For example let $(X, \calE)$ be a paved set that is stable for $(\cap f, \cup f)$; we will see in Proposition~2 that the family of $\calE$-measurable functions $\calM(X, \calE)$ is a convex cone of functions stable for $(\wedge f, \vee f)$. The following conditions are equivalent
\begin{enumerate}[i)]
\item $\delta(A \cap B) + \delta(A \cup B) = \delta(A) + \delta(B), \forall A, B \in \calE$,

\item $\int_X (f + g) d\delta = \int_X f d\delta + \int_X g d\delta, \forall f, g \in \calM^+(X, \calE)$,

\item $\int_X (f \wedge g) d\delta + \int_X (f \vee g) d\delta = \int_X f d\delta + \int_X g d\delta, \forall f,g \in \calM^+(X, \calE)$.
\end{enumerate}
If in i), ii) and iii) we change ``$=$'' to ``$\leqslant$'' or ``$\geqslant$'', we again obtain the equivalence of i), ii) and iii) (see~\cite{greco_preprint}).

\section{Properties of $\calE$-measurable functions.}

From the definition of $\calE$-measurability and Property ii) of Theorem~1, we obtain the following propositions on the properties of $\calE$-measurable functions.

{\it \textsc{Proposition 1}. Let $(X, \calE)$ be a paved set and let $\psi: [0, +\infty] \rightarrow [0, +\infty]$ be a continuous and increasing function. We have

\begin{enumerate}[i)]
\item if $\psi(0) = 0$, $f \in \calM^+(X, \calE)$, then $\psi \circ f \in \calM^+(X, \calE)$ (note we can dispense with the equality $\psi(0) = 0$ if $X \in \calE$);

\item if the sequence $\{ f_n \}_{n \in \N} \subset \calM^+(X, \calE)$ converges uniformly to $f$, then $f \in \calM^+(X, \calE)$.
\end{enumerate}}

{\it \textsc{Proposition 2}. Let $(X, \calE)$ be a paved set. Then

\begin{enumerate}[i)]
\item for $\calM^+(X, \calE)$ to be stable for $(\wedge f)$ (respectively, for $(\vee f)$) it is necessary and sufficient that $\calE$ be stable for $(\cap f)$ (respectively, for $(\cup f)$);

\item to have $f + g \in \calM^+(X, \calE)$ for any $f, g \in \calM^+(X, \calE)$, it is necessary and sufficient that $\calE$ be stable for $(\cap f, \cup f)$;

\item for $\calM^+(X, \calE)$ to be stable for $(\wedge d)$ (respectively, for $(\vee d)$) it is necessary and sufficient that $\calE$ be stable for $(\cap d)$ (respectively, for $(\cup d)$); in that case, for $f \in [0, +\infty]$ to be $\calE$-measurable it is necessary and sufficient that the set $\{ f \geqslant t \} \in \calE$ (respectively, $\{ f > t \} \in \calE$) for all $t \in (0, +\infty)$.
\end{enumerate}}

\textsc{Proof}. We are only going to show property ii). If $f + g$ is $\calE$ measurable for all pairs of functions $f, g \in \calM^+(X, \calE)$, we have that $\varphi_A + \varphi_B$ is $\calE$ measurable for all $A, B \in \calE$; it follows therefore that $A \cap B$ and $A \cup B \in \calE$, according to~ii) of Theorem~1; that is $\calE$ is stable for $(\cap f, \cup f)$. On the other hand, suppose that $\calE$ is stable for $(\cap f, \cup f)$, let us verify that $\calM^+(X, \calE)$ is stable for summation. If $g_1 = \sum_{i = 1}^n a_i \varphi_{A_i}$ and $g_2 = \sum_{i = 1}^n b_i \varphi_{B_i}$ belong to $\calS^+(X, \calE)$, the function $g_1 + g_2$ belongs to $\calS^+(X, \calE)$, because we have for all $t \in (0, +\infty)$ the equality
\[
\{ g_1 + g_2 > t \} = \cup_{(a, b) \in Z} (\{ g_1 \geqslant a \} \cap \{ g_2 \geqslant b \}),
\]
where $Z$ is the finite set
\[
\{ (g_1(x), g_2(x)) \in \R^2: x \in X \mbox{ and } g_1(x) + g_2(x) > t \}.
\]
If $f, g \in \calM^+(X, \calE)$ are bounded, there exist two sequences of bounded $\calE$-measurable functions $\{ g_n \}_{n \in \N}$, $\{ f_n \}_{n \in \N}$, which converge uniformly, respectively, to $g$ and $f$ (see Th.~1, property iii)); since $\{ g_n + f_n \}_{n \in \N} \subset \calS^+(X, \calE)$ converges uniformly to $g + f$, we therefore have that $f + g \in \calM^+(X, \calE)$, by property ii) of Proposition~1. Finally, taking into account the equality $\int_X (f + g) d\alpha = \lim_n \int_X (f \wedge n + g \wedge n) d\alpha$, true for any $f, g \in [0, +\infty]^X$ and for any monotonic set function $\alpha: \calP(X) \rightarrow [0, +\infty]$, it follows then that $f + g \in \calM^+(X, \calE)$, also if $f, g \in \calM^+(X, \calE)$ are not bounded. $\quad \quad \quad \blacksquare$

In particular cases, the class $\calM^+(X, \calE)$ has new properties. For example, if $X$ is a topological space in which every countable open cover of $X$ contains a finite open cover of $X$, any function $f \in [0, +\infty)^X$ that is measurable with respect to the paving of closed sets of $X$ (that is, $f \in [0, +\infty)^X$ is upper semi-continuous) is bounded; or if $X$ is a pseudocompact topological space (that is, every continuous function $g: X \rightarrow \R$ is bounded), any function $f \in [0, \infty)^X$ that is measurable with respect to the paving comprised of the sets of the form $\{ g = 0 \}$, where $g: X \rightarrow \R$ is continuous, is also bounded. This depends on the {\em semi-compactness} of the pavings in question (see~\cite{varadarayan1965}); recall that the paving $\calE$ on $X$ is said to be semi-compact (see~\cite{meyer1966}), if $\{ H_n \}_{n \in \N} \subset \calE$ and $\cap_{n = 1}^{\infty} H_n = \emptyset$ entails that there exists $n_0 \in \N$ such that $\cap_{n = 1}^{n_0} H_n = \emptyset$.

{\it \textsc{Theorem 2}. Let $X$ be a set endowed with a paving $\calE$ that is stable for $(\cap f)$. The following conditions are equivalent:

\begin{enumerate}[i)]
\item $\calE$ is semi-compact,
\item any function $f \in \calM^+(X, \calE) \cap \R^X$ is bounded,
\item any function $f \in \calM^+(X, \calE)$ has a maximum,
\item any sequence $\{ f_n \}_{n \in \N}$, decreasing to zero, converges uniformly to zero, if all $f_n$ are $\calE$-measurable,
\item for any sequence $\{ f_n \}_{n \in \N} \subset \calM^+(X, \calE)$, decreasing to zero, and for any monotonic $\alpha: \calP(X) \rightarrow [0, +\infty]$ such that $\int_X f_1 d\alpha < + \infty$, we have \linebreak $\lim_n \int_X f_n d\alpha = 0$.
\end{enumerate}}

\textsc{Proof}. It suffices to suitably adapt the proof of an analogous theorem of I.~Glicksberg~\cite{glicksberg1952}. i) $\implies$ iii) by property~ii) of Theorem~1, there exist sets $H_n \in \calE$ such that $\{ f \geqslant a_{n+1} \} \subset H_n \subset \{ f > a_n \}$, for all $n \in \N$, if $f \in \calM^+(X, \calE)$, $\{ a_n \}_{n \in \N} \subset (0, +\infty)$, $a_{n+1} > a_n$ for all $n \in \N$ and $\lim_n a_n = \sup \{ f(x): x \in X \}$ (we assume that $f \neq 0$); since $\cap_{n = 1}^m H_n \supset \{ f \geqslant a_{m+1} \} \neq \emptyset$ for all $m \in \N$, we deduce from the semi-compactness of $\calE$ that $\cap_{n = 1}^{\infty} H_n \neq \emptyset$; if $x_0 \in \cap_{n = 1}^{\infty} H_n$, we therefore have $f(x_0) \geqslant a_n$ for all $n \in \N$, that is $f(x_0) = \sup (f)$.

iii) $\implies$ ii): is an immediate consequence of iii).

ii) $\implies$ i): let $\{ H_n \}_{n \in \N} \subset \calE$ and $\cap_{n = 1}^{\infty} H_n = \emptyset$. The sets $G_n = \cap_{i = 1}^n H_n$ belong to $\calE$ and the function $\sum_{n = 1}^{\infty} n \varphi_{G_n} \in \calM^+(X, \calE) \cap \R^X$; by hypothesis ii) the function $\sum_{n = 1}^{\infty} n \varphi_{G_n}$ is bounded; then there exists a number $n_0 \in \N$ such that $G_{n_0} = \emptyset$, that is $\cap_{n = 1}^{n_0} H_n = \emptyset$.

i) $\implies$ iv): let $\{ f_n \}_{n \in \N} \subset \calM^+(X, \calE)$ be a sequence decreasing to zero. Let us show that $\{ f_n \}_n$ converges to zero uniformly, that is for any real number $a > 0$ there exists $n_0 \in \N$ such that $f_{n_0} \leqslant a$. We choose a sequence $\{ a_n \}_n \subset (0, +\infty)$ such that $a_n < a$ and $\lim_n a_n = a$. By the measurability of all $f_n$ and property ii) of Theorem~1, there exists for all $n \in \N$ a set $H_n \in \calE$ such that $\{ f_n \geqslant a \} \subset H_n \subset \{ f_n > a_n \}$. We therefore have $\cap_{n = 1}^{\infty} H_n = \emptyset$, because $\lim_n f_n = 0$; we then deduce that there exists $n_0 \in \N$ such that $\cap_{n = 1}^{n_0} H_n = \emptyset$, by hypothesis i). Therefore $\{ f_{n_0} \geqslant a \} \subset \cap_{n = 1}^{n_0} H_n = \emptyset$, that is the sequence $\{ f_n \}_{n \in \N}$ converges uniformly to zero.

iv) $\implies$ v): this is a consequence of the properties of the integral $\int_X \mbox{--- } d\alpha$ (see~\cite{greco1977}, Prop.~1, property~(4)).

v) $\implies$ i): it suffices to demonstrate that there exists $n_0 \in \N$ such that $H_{n_0} = \emptyset$, if $\{ H_n \}_{n \in \N} \subset \calE$ with $\cap_{n = 1}^{\infty} H_n = \emptyset$ and $H_{n+1} \subset H_n$ for all $n \in \N$. We reason by contradiction; if $H_n \neq \emptyset$ for all $n \in \N$, the monotone set function $\alpha_{\calF}: \calP(X) \rightarrow \{ 0, 1 \}$ that is equal to 1 only for those sets of the filter $\calF$ generated by the $H_n$, is such that $\alpha_{\calF} = 1$ for all $n \in \N$; which is contrary to hypothesis v), since the sequence $\{ \varphi_{H_n} \}_{n \in \N}$ decreases to zero. $\quad \quad \quad \blacksquare$

When the paving $\calE$ on $X$ is an algebra of sets we have the following characterisation of $\calE$-measurability.

{\it \textsc{Theorem 3}. Let $(X, \calA)$ be a paved set, $\calA$ an algebra of sets, $f \in \overline{\R}^x$. The following conditions are equivalent:

\begin{enumerate}[i)]
\item $f$ is $\calA$-measurable,
\item for any ultrafilter $\calF$ of the algebra $\calA$ there exists the limit $\lim_{\calF} f \in \overline{\R}$,
\item for any real number $a > 0$, there exists a finite partition $\{ H_i \}_{i = 1}^n \subset \calA$ of $X$ such that for all $i$ we have
\[
f(x) \leqslant a + f(y) \quad \mbox{if } x, y \in H_i
\]
or
\[
f(x) \geqslant \frac{1}{a} \quad \mbox{if } x \in H_i
\]
or
\[
f(x) \leqslant -\frac{1}{a} \quad \mbox{if } x \in H_i.
\]
\end{enumerate}}

\textsc{Proof}. i) $\implies$ ii): we can limit ourselves to verifying that for all $f \in [0, +\infty]^X \cap \calM(X, \calE)$ there exists the limit $\lim_{\calF} f \in \overline{\R}$. Let $\calF' \subset \calP(X)$ be the filter of $\calP(X)$ generated by the ultrafilter $\calF$ of $\calA$, and let $\alpha_{\calF'}$ and $\beta_{\calF'}$ be the monotonic set functions defined in the introduction. For these set functions we have $\int_X f d\alpha_{\calF'} = \liminf_{\calF} f$ and $\int_X f d\beta_{\calF'} = \limsup_{\calF} f$; since $\alpha_{\calF'}(H) = \beta_{\calF'}(H)$ for any $H \in \calA$, it follows that $\limsup_{\calF} f = \liminf_{\calF} f$, by the definition of $\calA$-measurability; this proves that $\lim_{\calF} f$ exists and belongs to $\overline{\R}$.

ii) $\implies$ iii): according to ii), for any ultrafilter $\calF$ of the algebra $\calA$ there exists a set $H_{\calF} \in \calF$\footnote{Translator's note: The original has $H_{\calF} \in \calA$.} such that $f(x) - f(y) \leqslant a$ or $f(x) \geqslant 1/a$ or else $f(x) \leqslant - 1/a$ for $x, y \in H_{\calF}$. The family $\{ X - H_{\calF}: \calF \mbox{ an ultrafilter of } \calA \}$ cannot be contained in any ultrafilter of $\calA$; there exists therefore a finite number $\calF_1, \calF_2, \ldots, \calF_n$ of ultrafilters of $\calA$ such that
\[
\cap_{i = 1}^n (X - H_{\calF_i}) = \emptyset ;
\]
the proof of the implication ii) $\implies$ iii) is thus complete.

iii) $\implies$ i): it suffices to verify this implication when $f \in [0, +\infty]^X$ is bounded. In that case the construction of a sequence $\{ g_n \}_n \subset \calM^+(X, \calA)$ such that $\{ g_n \}_n$ converges uniformly to the function $f$ is immediate; thus we can deduce that $f$ is $\calA$-measurable, by the properties of $\calA$-measurable functions. $\blacksquare$

\textsc{Remark~2}. Let $X = \R$ and let $\calA_{\R}$ be the algebra generated by the bounded and unbounded intervals of $\R$. For all $x \in \R$, let $\calF_x$, $\calF_{x^+}$, $\calF_{x^-}$, $\calF_{+\infty}$, $\calF_{-\infty}$ be the ultrafilters of $\calA_{\R}$ generated, respectively, by the set $\{ x \}$, by the intervals $(x, x + a)$, where $a \in (0, +\infty)$, by the intervals $(x - a, x)$, where $a \in (0, +\infty)$, by the intervals $(a, +\infty)$, where $a \in \R$, and by the intervals $(-\infty, a)$, where $a \in \R$; all the ultrafilters of $\calA_{\R}$ are of these types. According to Theorem~3, it follows that a function $f \in \overline{\R}^X$ is $\calA_{\R}$-measurable if and only if the following limits exist: $\lim_{y \rightarrow x^+} f(y)$, $\lim_{y \rightarrow x^-} f(y)$, $\lim_{y \rightarrow +\infty} f(y)$, and $\lim_{y \rightarrow -\infty} f(y)$; that is $f$ is devoid of oscillatory discontinuities. This example can easily be adapted when $\calA$ is the algebra generated by the intervals that are semi-open on the right (or on the left), and when $X = \R^n$ and $\calA$ is the algebra generated by the $n$-dimensional intervals.

\textsc{Remark~3}. Theorem~3 can be presented in a more suggestive manner. If we denote by $[\calA]$ the Stone space (the space of ultrafilters of $\calA$)(see~\cite{sikorski1969}), associated with the algebra $\calA$, for a function $f \in \overline{\R}^X$ to be $\calA$-measurable it is necessary and sufficient that there exists a continuous function $\overline{f}: [\calA] \rightarrow \overline{\R}$ such that $\overline{f}(\calF_x) = f(x)$ for any ultrafilter $\calF_x = \{ H \in \calA: x \in H \}$, where $x \in X$. The class of bounded and $\calA$-measurable functions is therefore an algebra of functions, isomorphic to the algebra of continuous functions from $[\calA]$ to $\R$.

\textsc{Remark~4}. In the particular case where the algebra $\calA$ is a $\sigma$-algebra, we have that $f + g \in \calM(X, \calA)$ if $f + g$ is defined over $X$ and $f, g \in \calM(X, \calE)$. This property is not true if $\calA$ is only an algebra. Indeed if we consider the algebra $\calA \subset \calP(\N)$ of finite and cofinite subsets of $\N$, for a function $f: \N \rightarrow \overline{\R}$ to be $\calA$-measurable it is necessary and sufficient that there exists the limit $\lim_n f(n)$; in that case we can therefore easily choose two functions $f, g: \N \rightarrow \overline{\R}$ such that $f + g$ is defined, but the function $f + g$ is not $\calA$-measurable.

\textsc{Remark~5}. Let $\mu: \calP(X) \rightarrow [0, +\infty]$ be a monotonic or non-monotonic set function with $\mu(\emptyset) = 0$. We know that the family of sets that are $\mu$-measurable in the sense of Carath\'eodory is a set algebra $\calA_{\mu} \subset \calP(X)$. Let $\delta: \calA_{\mu} \rightarrow [0, +\infty]$ be the monotone set function defined by $\delta(H) = \mu(H)$ for all $H \in \calA_{\mu}$; a function $f$ is said to be $\mu$-measurable if $f$ is $\calA_{\mu}$-measurable. For these $\mu$-measurable functions we have linearity, that is $\int (f + g) d\delta = \int_X f d\delta + \int_X g d\delta$ (see Remark~1 for the definition of the integral $\int_X \mbox{--- } d\delta$) for all $f, g \in \calM^+(X, \calA_{\mu})$. Furthermore we have that the sets $\{ t \in \R: \{ f > t \} \notin \calA_{\mu} \}$ and $\{ t \in \R: \{ f \geqslant t \} \notin \calA_{\mu} \}$ are countable, if $f \in \calM(X, \calA_{\mu})$ and $\int_X | f | d\delta < + \infty$. The class of functions $f \in \R^X$ that are $\mu$-measurable, such that $\int_X | f | d\delta < + \infty$, is a Riesz space (see~\cite{greco_preprint}).

\section{Measurable functions with respect to two pavings.}

In this paragraph we consider two pavings $\calK, \calU$ on $X$, both stable for $(\cap f, \cup f)$; we define on $\calP(X)$ the relation ``$\ll$'' in this manner: ``$A \ll B$ if and only if there exists a pair $(K, U) \in \calK \times \calU$ such that $A \subset K \subset U \subset B$''. Let $D$ be the set of dyadic numbers in the interval $[0, 1]$; let $\{ H_t \}_{t \in D} \subset \calP(X)$ be a family of sets such that $H_0 = X$ and $H_t \gg H_s$ for all $t, s \in D$ with $t < s$; then the function $f: X \rightarrow [0, 1]$, defined by $f(x) = \sup \{ t \in D: x \in H_t \}$, is either $\calK$-measurable or $\calU$-measurable. Indeed for any $a, b \in (0, +\infty)$ with $a > b$, let $t$ and $r$ be two numbers $\in D$ such that $a > t > r > b$, it therefore follows that $\{ f \geqslant a \} \subset H_t \ll H_r \subset \{ f > b \}$; since there exists a pair $(K, U) \in \calK \times \calU$ such that $H_t \subset K \subset U \subset H_r$, we finally have that $\{ f \geqslant a \} \subset K \subset U \subset \{ f > b \}$, that is $f$ is $\calK$-measurable and $\calU$-measurable.

We say that the pair of pavings $(\calK, \calU)$ satisfies {\em property} ($N$) if for any pair $(K, U) \in \calK \times \calU$ with $K \subset U$, there exists a pair $(K', U') \in \calK \times \calU$ such that $K \subset U' \subset K' \subset U$.

{\it \textsc{Lemma~3}. Let $k, u \in [0, +\infty]^X$, and let $(\calK, \calU)$ be a pair of pavings on $X$ that satisfy property ($N$). If $k \leqslant u$ and if $\{ k \geqslant t \} \in \calK$ and $\{ u > t \} \in \calU$ for any $t \in (0, +\infty)$, then there exists a function $f \in \calM^+(X, \calK) \cap \calM^+(X, \calU)$ such that $k \leqslant f \leqslant u$.}

\textsc{Proof}. By property~i) of Proposition~1 we can suppose that $u \leq 1$. For any $t \in D$, where $D$ is the set of dyadic numbers in the interval $[0, 1]$, we set $K_t = \{ k \geqslant t \}$ and $U_t = \{ u > t \}$. We construct a family of sets $\{ F_t \}_{t \in D} \subset \calP(X)$, enjoying the following properties\footnote{Translator's note: In (2), the original has $G_t$ in place of $U_t$, and $H_r$ in place of $K_r$.}:
\begin{eqnarray}
F_0 &=& X, \quad F_1 = \emptyset, \\
F_r &\ll& F_t, \quad F_r \ll U_t, \quad K_r \ll F_t, \quad \mbox{if } t, r \in D \mbox{ and } t < r ;
\end{eqnarray}
in fact if we suppose sets $F_t$ satisfying~(1) and~(2) have been defined for all $t = i / 2^n$, where $0 \leqslant i \leqslant 2^n$, then we can choose a set $F_{t_0}$ for $t_0 = (2i + 1)/2^{n + 1}$, where $0 \leqslant i \leqslant 2^n - 1$, satisfying the relations: $F_{i / 2^n} \gg F_{t_0}$, $F_{t_0} \gg F_{(i + 1) / 2^n}$, $U_{i / 2^n} \gg F_{t_0}$, $F_{t_0} \gg K_{(i + 1)/2^n}$, according to the relations (true by hypothesis or by construction):
\[
F_{i/2^n} \gg F_{(i+1) / 2^n}, \; F_{i / 2^n} \gg K_{(i + 1) / 2^n}, \; U_{i / 2^n} \gg F_{(i = 1) / 2^n}, \; U_{i / 2^n} \gg K_{(i + 1) / 2^n}.
\]
The function $f: X \rightarrow [0, +\infty]$, defined by $f(x) = \sup \{ t \in D: x \in F_t \}$, is thus $\calK$-measurable and $\calU$-measurable and $k \leqslant f \leqslant u$. $\quad \quad \quad \blacksquare$

By Lemma~3 and property~iii) of Proposition~2 we deduce the following theorems:

{\it \textsc{Theorem~4}. Let $\calK$ and $\calU$ be two pavings on $X$, stable for $(\cap f, \cup f)$. The following properties are equivalent:
\begin{enumerate}[i)]
\item the pair $(\calK, \calU)$ satisfies property ($N$),

\item for any pair $(K, U) \in \calK \times \calU$ there exists a function $f \in \calM^+(X, \calK) \cap \calM^+(X, \calU)$ such that $\varphi_K \leqslant f \leqslant \varphi_U$, if $K \subset U$. $\quad \quad \quad \blacksquare$
\end{enumerate}}

{\it \textsc{Theorem~5}. Let $\calK$ be a paving stable for $(\cap d, \cup f)$ and let $\calU$ be a paving stable for $(\cap f, \cup d)$. The following properties are equivalent:
\begin{enumerate}[i)]
\item the pair $(\calK, \calU)$ satisfies property ($N$),

\item for any $(k, u) \in \calM^+(X, \calK) \times \calM^+(X, \calU)$ there exists a function \linebreak $f \in \calM^+(X, \calK) \cap \calM^+(X, \calU)$ such that $k \leqslant f \leqslant u$, if $k \leqslant u$. $\quad \quad \quad \blacksquare$
\end{enumerate}}

Let $X$ be a topological space and let $\calK$ (respectively, $\calU$) be a paving on $X$ comprised of closed sets (respectively, open sets); the pair of pavings $(\calK, \calU)$ satisfies property~($N$) if and only if $X$ is normal. In that case Theorems~4 and~5 are well known theorems on normal spaces (see~\cite{tong1952}).

\end{document}